\documentclass[11pt, a4paper, reqno]{amsart}
\usepackage{preambolo}

\begin{document}


\setcounter{secnumdepth}{3}

\setcounter{tocdepth}{2}

\title{\textbf{Kodaira dimension of $\su(\lowercase{m})$-structures}
}

\author[Lorenzo Sillari]{Lorenzo Sillari}

\address{Lorenzo Sillari: Dipartimento di Scienze Matematiche, Fisiche e Informatiche, Unità di Matematica e Informatica, Università degli Studi di Parma, Parco Area delle Scienze 53/A, 43124, Parma, Italy} \email{lorenzo.sillari@unipr.it}

\author[Adriano Tomassini]{Adriano Tomassini}

\address{Adriano Tomassini: Dipartimento di Scienze Matematiche, Fisiche e Informatiche, Unità di Matematica e Informatica, Università degli Studi di Parma, Parco Area delle Scienze 53/A, 43124, Parma, Italy}
\email{adriano.tomassini@unipr.it}

\maketitle

\begin{abstract} 
\noindent \textsc{Abstract}. We study the Kodaira dimension of almost complex manifolds admitting an $\su (m)$-structure. We introduce the notion of almost complex structure of splitting type and of associated $\su(m)$-structure. When the latter is pseudoholomorphic, we provide two constructions that allow to obtain non-invariant almost complex structures with Kodaira dimension $0$, resp.\ with Kodaira dimension $-\infty$. Our results apply, in particular, to complex structures of splitting type and to several almost complex manifolds already well-studied in the literature. 
\end{abstract}

\blfootnote{  \hspace{-0.55cm} 
{\scriptsize 2020 \textit{Mathematics Subject Classification}. Primary: 32Q60, 53C10; Secondary: 32L05, 53C15. 
\vspace{.1cm}

\noindent The authors are partially supported by GNSAGA of INdAM. The second author is partially supported by the Project PRIN 2022 “Real and Complex Manifolds: Geometry and Holomorphic Dynamics” (code 2022AP8HZ9)}}

\section{Introduction}\label{sec:intro}

The Kodaira dimension is a classical complex invariant that measures the growth rate of the space of sections of the pluricanonical bundle of a complex manifold. It becomes especially useful in low dimension, where it plays a significant role in Enriques--Kodaira's classification of complex surfaces.

Recently, Chen and Zhang gave a definition of \textit{Kodaira dimension of almost complex manifolds} that, in the integrable case, reduces to the usual one, see \cite{CZ23} and \cite{CZ24}. They study the basic properties of the Kodaira dimension (in particular on $4$-manifolds) and of the spaces of pseudoholomorphic sections of \textit{bundle almost complex structures} \cite{BT96}. Furthermore, they provide several explicit examples and computations and they ask for other explicit computations of Kodaira dimension, possibly on manifolds with a certain structure (e.g., fibrations). So far, the available explicit computations are:
\begin{itemize}
    \item on the $4$-torus and on $4$-manifolds that admit a pseudoholomorphic elliptic fibration over Riemann surfaces of genus $g \ge 2$ \cite{CZ23};

    \item  on the Kodaira--Thurston $4$-manifold, see \cite{Cat22} and \cite{CZ23}, or on its higher-dimensional generalization \cite{HP24};

    \item on $4$-dimensional solvmanifolds without complex structures \cite{CNT21};

    \item on the completely solvable Nakamura manifold \cite{CNT20};

    \item on the Lie groups $\su(2) \times \su(2)$, $\mathrm{SO} (4)$, $\mathrm{Sp}(2)$, $\su(3)$, and on Lie groups of the form $G \times G$, where $G$ is a connected Lie group \cite{CNT24};
    
    \item on the $6$-sphere \cite{CZ23}.
\end{itemize}

All of the mentioned computations, with the exception of the $4$-torus, are performed for left-invariant structures on Lie groups or for invariant structures on compact quotients of Lie groups. Even though determining the Kodaira dimension of an invariant structure is already a hard task, finding that of a non-invariant one usually requires solving a difficult PDE whose coefficients are smooth functions (or showing that such a PDE has no solution).
\vspace{.2cm}

Since the almost complex Kodaira dimension can be computed in terms of sections of the canonical bundle, it is natural to relate it to $\su(m)$-structures, see, e.g., \cite{Bry87} or \cite{Sal89}. An almost complex $2m$-manifold admits an associated $\su(m)$-structure when its canonical bundle is trivial as a complex (not necessarily pseudoholomorphic) line bundle. Lie groups and their compact quotients by lattices are parallelizable, hence they always admit an $\su(m)$-structure, and it becomes natural to exploit its existence to help in the computations of the Kodaira dimension.
\vspace{.2cm}

In this paper, we focus on the Kodaira dimension of $\su(m)$-structures. We have two main goals. The first one, is to study general properties of the Kodaira dimension of such structures. We are able to relate the vanishing of the Kodaira dimension to the triviality of tensor powers of the canonical bundle, see Section \ref{sec:prel}. Our results apply to parallelizable manifolds and to compact quotients of Lie groups, and allow to recover results from \cite{AT24} and \cite{CNT24} as special cases.

The second goal is to partially answer the problem posed by Chen and Zhang by providing two constructions that allow to build \textit{non-invariant} almost complex structures for which we can compute the Kodaira dimension. To do so, in Section \ref{sec:prel} we introduce \textit{almost complex structures of splitting type} on Lie groups and their compact quotients, see Definition \ref{def:splitting}, a natural non-integrable broad generalization of complex structures of splitting type on solvmanifolds, see \cite{AOUV17} and \cite{Kas13}. Compact quotients of these Lie groups have the structure of a pseudoholomorphic fibration over a complex torus, and they inherit a standard $\su(m)$-structure compatible with the fibration. Section \ref{sec:fibrations} is where we prove our main results: when the standard $\su(m)$-structure has pseudoholomorphically trivial canonical bundle, we use the fibration to build many non-invariant almost complex structures with Kodaira dimension $0$ or $-\infty$. More precisely, we have the following.

\begin{mythm}{\ref{thm:family:0}}
    Let $G = \C^k \ltimes_\eta H$ be a $2m$-dimensional Lie group endowed with an almost complex structure of splitting type whose standard $\su (m)$-structure is pseudoholomorphic, and let $M = \Gamma \backslash G$ be its compact quotient. Then there exists a family of non-invariant almost complex structures $J_F$ on $M$ with $\kappa_{J_F} =0$, parametrized by $F \in (\mathcal{F})^{m-k}$. 
\end{mythm}

\begin{mythm}{\ref{thm:family:infty}}
    Let $G = \C^k \ltimes_\eta H$ be a $2m$-dimensional Lie group endowed with an almost complex structure of splitting type whose standard $\su (m)$-structure is pseudoholomorphic, and let $M = \Gamma \backslash G$ be its compact quotient. Then there exists a family of non-invariant almost complex structures $J_F$ on $M$ with $\kappa_{J_F} = - \infty$, parametrized by $F \in C^\infty(M)^k$. 
\end{mythm}

In the statement of the theorems, we denoted by $\Gamma$ a lattice in $G$ and by $\mathcal{F}$ a suitable subspace of functions in $C^\infty(M)$. The assumptions of the theorems hold, in particular, for \textit{complex} structure of splitting type, see Remark \ref{rem:integrable}.

We point out that all the examples that we mentioned, on which the Kodaira dimension has been computed (except for the $6$-sphere), can be seen as Lie groups endowed with an almost complex structure of splitting type or as their compact quotients. The same is true for other relevant manifolds present in the literature, whose Kodaira dimension has not been computed before. In Section \ref{sec:examples}, we show some application of our results to those examples. In particular, our construction can be performed on $2m$-dimensional tori, on the Iwasawa manifold, on the holomorphically parallelizable Nakamura manifold, on the $6$-dimensional solvmanifolds without invariant complex structures built by Fern\'andez, de Le\'on and Saralegui \cite{FLS96}, and on the $8$-dimensional solvmanifold built by Sferruzza and the second author \cite{SfeT24}.

\section{Preliminaries}\label{sec:prel}

\subsection*{Notation.} Given a basis of complex forms $\phi^j$, $j = 1, \ldots, m$, we abbreviate $\phi^j \wedge \phi^k$ to $\phi^{jk}$ and $\bar \phi^j$ to $\phi^{\bar j}$.
\vspace{.2cm}

In this section, we recall some preliminaries on almost complex structures and their Kodaira dimension, on $\su(m)$-structures, and on the canonical bundle of compact quotients of Lie groups.

\subsection{Almost complex Kodaira dimension.}\label{sec:kodim}

Let $(M,J)$ be a compact almost complex $2m$-manifold. Denote by $\Lambda^{p,q}$ the bundle of complex $(p,q)$-forms and by $A^{p,q}$ the space of its smooth sections, that is, the space of complex $(p,q)$-forms.

The \textit{canonical bundle} of $(M,J)$ is the bundle $K_J \coloneqq \Lambda^{m,0}$. A $(p,q)$-form $\alpha \in A^{p,q}$ is called \textit{pseudoholomorphic} if $\bar \partial \alpha=0$. In particular, a pseudoholomorphic section of the canonical bundle is a $\bar \partial$-closed complex $(m,0)$-form $\Omega$. If $J$ is integrable, pseudoholomorphic forms are actually holomorphic forms. Note that, in the complex case, the term \textit{canonical bundle} usually refers to the bundle of \textit{holomorphic $(m,0)$-form}, see \cite{Huy05}. Even though in practice we will often consider pseudoholomorphic $(m,0)$-forms, we do not require $\bar \partial$-closedness in the definition of $K_J$.

Denote by $P_l$ the space of pseudoholomorphic sections of $K_J^{\otimes l}$, for $l \ge 1$. The space $P_1$ is just the space of pseudoholomorphic $(m,0)$-forms. Chen and Zhang recently introduced the notion of \emph{Kodaira dimension of an almost complex manifold}, see \cite{CZ23} and \cite{CZ24}. The almost complex Kodaira dimension is the number
    \[
    \kappa_J := 
    \begin{cases}
        - \infty &\text{if $P_l = \{ 0 \}$ for all $l \ge 1$,}\\
        \limsup\limits_{l \to \infty} \frac{\log \dim_\C P_l}{\log l} &\text{otherwise},
    \end{cases}
    \]
and it generalizes the usual notion of Kodaira dimension of a complex manifold. 

\begin{remark}
In addition to the definition given above, Chen and Zhang gave another, more algebraic, definition of Kodaira dimension of almost complex manifolds. In the complex case, the two definitions are equivalent to the usual definition of Kodaira dimension. In the non-integrable case, it is not known if the two definitions are equivalent. Due to the analytical approach of this paper, we deal with the definition given in terms of sections of tensor powers of the canonical bundle that we recalled above.
\end{remark}

\subsection{\texorpdfstring{$\su (m)$}{}-structures.}\label{sec:sum}

Let $M$ be a compact $2m$-manifold. We say that $M$ admits an $\su (m)$-structure if the structure group of its tangent bundle can be reduced to $\su (m)$. The choice of a reduction corresponds to the choice of an $\su (m)$-structure, which determines an almost complex structure $J$ whose canonical bundle $K_J$ is trivial, and the choice of a trivialization $\Omega$ of $K_J$, that is, a never-vanishing $(m,0)$-form. This also provides a preferred volume form on $M$, given by $\Omega \wedge \bar \Omega$. If, in addition, we can take $\Omega$ to be pseudoholomorphic, then $K_J$ is pseudoholomorphically trivial. By Proposition 3.6 in \cite{CNT24}, if the canonical bundle is holomorphically trivial, then the Kodaira dimension of $J$ vanishes. We show that this is also the case when a tensor power of the canonical bundle is holomorphically trivial.

\begin{lemma}\label{lemma:kodaira:holomorphic}
    Let $(M,J)$ be a compact almost complex manifold admitting an $\su (m)$-structure. If the canonical bundle of $J$ or one of its tensor powers is pseudoholomorphically trivial, then $\kappa_J = 0$.
\end{lemma}
\begin{proof}
    Let $\Omega$ be a (not necessarily pseudoholomorphic) trivialization of $K_J$ given by an $\su(m)$-structure, and let $\gamma$ be the $(0,1)$-form such that $\bar \partial \Omega = \gamma \wedge \Omega$. The trivialization $\Omega$ is pseudoholomorphic if and only if $\gamma$ vanishes. Our claim for the case $\gamma=0$ is already proved in Proposition 3.6 in \cite{CNT24} for the case of parallelizable manifolds. The proof for $\su(m)$-manifolds is very similar, therefore we omit it, and we can suppose $\gamma \neq 0$. A section $f \Omega^{\otimes l}$ of $K_J^{\otimes l}$ is pseudoholomorphic if and only if $f$ is a solution of the equation
    \[
    D_l (f) \coloneqq \bar \partial f + l f \gamma = 0,
    \]
    cf.\ \cite{CNT24} or \cite{CZ23}. By assumption, there exists a minimum integer $l_0$ such that $K_J^{\otimes l_0}$ is pseudoholomorphically trivial. We show that
    \[
    \dim_\C P_l =
    \begin{cases}
        1, & \text{if $l$ is a multiple of $l_0$},\\
        0, &\text{otherwise},
    \end{cases}
    \]
    which readily implies that $\kappa_J =0$.

    Let $F \Omega^{\otimes l_0}$ be a trivialization for $K_J^{\otimes l_0}$, where $F$ is a never-vanishing smooth function. Then $F$ solves the equation $D_{l_0} (F) =0$. Let $g \Omega^{\otimes l_0}$ be another pseudoholomorphic section of $K_J^{\otimes l_0}$, where $g$ is not necessarily never-vanishing. Then 
    \[
    0 = \bar \partial (g \Omega^{\otimes l_0}) = \bar \partial (\frac{g}{F} F \Omega^{\otimes l_0}) = \bar \partial(\frac{g}{F}) \wedge ( F \Omega^{\otimes l_0}).
    \]
    Since $F \Omega^{\otimes l_0}$ is never-vanishing, this implies that $\bar \partial (g/F) =0$ and, by compactness of $M$, that $g/F$ is constant. Hence, we have that $\dim_\C P_{l_0} = 1$. Similarly, if $l = p l_0$, then we have that $D_l (F^p) =0$, that $F^p \Omega^{\otimes l}$ is never-vanishing and pseudoholomorphic, and that every other pseudoholomorphic section of $K_J^{\otimes l}$ is a constant multiple of $F^p \Omega^{\otimes l}$, which gives $\dim_\C P_l = 1$.
    
    Now suppose, by contradiction, that for $p < l_0$ there exists a function $g$ such that $ g \Omega^{\otimes p}$ is pseudoholomorphic. Then $g^{l_0} \Omega^{\otimes p l_0}$ and $ F^p \Omega^{\otimes pl_0}$ are both pseudoholomorphic sections of $K_J^{\otimes p l_0}$. By the argument of the first part of the proof, we have that $g^{l_0}$ is a constant multiple of $F^p$, which is never-vanishing. In particular, $g^{l_0}$ is never-vanishing, and so is $g$. This contradicts the minimality assumption on $l_0$ and proves that $K_J^{\otimes p}$ has no pseudoholomorphic sections for $p < l_0$, giving $\dim_\C P_p = 0$. 
    
    To conclude the proof, let $p$ be any integer that is not a multiple of $l_0$, and let $q$ be the integer such that $q l_0 < p < (q+1) l_0$. Suppose that $g \Omega^{\otimes p}$ is a pseudoholomorphic section of $K_J^{\otimes p}$. Then, we have that
    \begin{align*}
    0 &= \bar \partial ( g \Omega^{\otimes p} ) = \bar \partial ( \frac{g}{F^q} \Omega^{\otimes (p - q l_0)} \otimes F^q \Omega^{\otimes q l_0} ) = \bar \partial ( \frac{g}{F^q} \Omega^{\otimes (p - q l_0)}) \otimes (F^q \Omega^{\otimes ql_0})\\
    & = (\bar \partial (\frac{g}{F^q}) + (p-ql_0) \frac{g}{F^q} \gamma) \wedge ( F^q \Omega^{\otimes p}),
    \end{align*}
    which is equivalent to $D_{p-ql_0} ( g/ F^q) =0$. Since $p - q l_0 < l_0$, there are no non-zero solutions to such equation, and we get $\dim_\C P_p = 0$.
\end{proof}

\subsection{Compact quotients of Lie groups with trivial canonical bundle.}

Let $G$ be a connected, simply connected Lie group admitting a co-compact lattice $\Gamma$, and denote by $M \coloneqq \Gamma \backslash G$ their compact quotient. If $G$ is nilpotent, resp.\ solvable, we say that $M$ is a \textit{nilmanifold}, resp.\ \textit{solvmanifold}. If $G$ is completely solvable, we say that $M$ is a solvmanifold \textit{of completely solvable type}. 

Since $G$ admits a lattice, it must be unimodular. This implies that the quotient $M$ is parallelizable. In particular, it admits invariant almost complex structures and invariant $\su(m)$-structures induced by structures on $\g$, the Lie algebra of $G$. In addition, it will also admit genuine almost complex or $\su (m)$-structures that are not invariant.

Let $J$ be an invariant almost complex structure on $M$, and let $\phi^j$, $j =1, \ldots, m$, be a basis of invariant $(1,0)$-forms for $J$. Then, the invariant $(m,0)$-form $\phi^{1 \ldots m}$ induces a trivialization of the canonical bundle $K_J$, hence a canonical invariant $\su (m)$-structure on $M$. By choosing a different basis of $(1,0)$-forms, the trivialization changes by a constant. 

As a direct application of Lemma \ref{lemma:kodaira:holomorphic}, we have the following.

\begin{proposition}\label{prop:kodaira:quotient}
    Let $M = \Gamma \backslash G$ be a compact quotient of a Lie group and let $J$ be an invariant almost complex structure on $M$. If the canonical bundle of $J$ or one of its tensor powers is pseudoholomorphically trivial, then $\kappa_J = 0$.
\end{proposition}

We recall that if $(M,J)$ is a nilmanifold endowed with an invariant \textit{complex} structure, then its canonical bundle is always holomorphically trivial \cite{BDV09}, so that $\kappa_J = 0$. Moreover, the trivialization is necessarily provided by an invariant form. In the solvable case, there are examples of invariant \textit{complex} structures whose canonical bundle is holomorphically trivial, but the trivialization cannot be invariant, see \cite{AT24} and \cite{Tol24}. We point out that this difference amounts to understanding if a given invariant complex structure admits or not a compatible invariant $\su (m)$-structure with holomorphically trivial canonical bundle.

Finally, we give the almost complex version of Corollary 4.2 in \cite{AT24}.

\begin{lemma}\label{lemma:invariant:trivialization}
    Let $M = \Gamma \backslash G$ be a compact quotient of a Lie group and let $J$ be an invariant complex structure on $M$. Then pseudoholomorphic trivializations of tensor powers of its canonical bundle are either all invariant or all non-invariant (if they exist).
\end{lemma}
\begin{proof}
    Let $\Omega$ be a (not necessarily pseudoholomorphic) trivialization of $K_J$. Suppose that the bundle $K_J^{\otimes l}$ admits a never-vanishing pseudoholomorphic section $F \Omega^{\otimes l}$, $F \in C^\infty (M)$, and let $\hat{F} \Omega^{\otimes l}$, $\hat{F} \in C^\infty (M)$, be another pseudoholomorphic section of $K_J^{\otimes l}$. Since $F \Omega^{\otimes l}$ trivializes $K_J^{\otimes l}$, there exists a function $g$ such that $\hat{F} \Omega^{\otimes l} = g F \Omega^{\otimes l}$. Then, we have that
    \[
    0 = \bar \partial ( \hat{F} \Omega^{\otimes l} ) = \bar \partial (g F \Omega^{\otimes l}) = F \bar \partial g \wedge \Omega^{\otimes l} + g \bar \partial ( F \Omega^{\otimes l} ) =  F \bar \partial g \wedge \Omega^{\otimes l},
    \]
    where in the last equality we used that $F \Omega^{\otimes l}$ is pseudoholomorphic. Since $F$ is never-vanishing and $\Omega^{\otimes l}$ induces an isomorphism between $A^{0,1}$ and $A^{0,1} \wedge K_J^{\otimes l}$, we conclude that $\bar \partial g =0$ and, by compactness of $M$, that $g$ is constant. In particular, every pseudoholomorphic section of $K_J^{\otimes l}$ is a constant multiple of $F \Omega^{\otimes l}$. The trivialization $F \Omega^{\otimes l}$ is invariant if and only if $F$ is constant, if and only if $g$ is constant. Similarly, $F$ is non-constant if and only if $g$ is. In this case every trivialization is non-invariant.
\end{proof}

\subsection{Almost complex and \texorpdfstring{$\su (m)$}{}-structures of splitting type.} 

We introduce a special class of Lie groups that split as the semi-direct product of $\C^k$ and a Lie group endowed with a left-invariant almost complex structure, with the action of $\C^k$ compatible with the almost complex structure.

Let $G$ be a connected, simply connected, $2m$-dimensional Lie group.

\begin{definition}\label{def:splitting}
    Let $G$ be a Lie group that can be written as a semi-direct product $G = \C^k \ltimes_\eta H$, $k \ge 1$, such that: 
    \begin{itemize}
        \item [(i)] $H$ is a connected, simply connected, $2(m-k)$-dimensional Lie group endowed with a left-invariant almost complex structure $J_H$;

        \item [(ii)] $G$ admits a co-compact lattice $\Gamma$ that splits as $\Gamma = \Gamma_{\C^k} \ltimes_\eta \Gamma_H$ and such that, for all $z \in \C^k$, we have $\eta(z) \in \Aut (\Gamma_H)$.
    \end{itemize} 
Then $G$ admits a standard left-invariant almost complex structure built as follows. Let $dz^j$, $j = 1, \ldots, k$, be the standard basis of $(1,0)$-forms for $\C^k$, let $\psi^j$, $j = k+1, \ldots, m$, be a basis of left-invariant $(1,0)$-forms on $H$ for the almost complex structure $J_H$ and set $\phi^j = \eta^{-1} \psi^j$, $j = k+1, \ldots, m$.

The \textit{almost complex structure of splitting type} induced on $G$ by $J_H$ is the one defined by the co-frame of $(1,0)$-forms 
\begin{equation}\label{eq:std:basis}
\{ dz^1, \ldots, dz^k, \phi^{k+1}, \ldots,\phi^{m} \}.
\end{equation}
Note that the construction depends only on the choice of co-frame of $(1,0)$-forms for the almost complex structure $J_H$. In this case, we say that $G$ is an \textit{almost complex Lie group of splitting type}.
\end{definition}
Every almost complex structure of splitting type naturally defines an $\su (m)$-structure given by the left-invariant $(m,0)$-form
\[
\Omega \coloneqq dz^{1 \ldots k} \wedge \phi^{k+1 \ldots m},
\]
that we call the \textit{standard left-invariant $\su (m)$-structure} associated to $J$. In general, the form $\Omega$ will not be pseudoholomorphic since
\[
\bar \partial \Omega = (-1)^k dz^{1 \ldots k} \wedge  \bar \partial ( \phi^{k+1 \ldots m} ).
\]
If $\bar \partial \Omega =0$, we say that the standard left-invariant $\su(m)$-structure on $G$ is \textit{pseudoholomorphic}. Note that pseudoholomorphic $\su (m)$-structures have already been considered, under different names, by de Bartolomeis and the second author in \cite{BT13} and \cite{BT06b}.

\begin{remark}\label{rem:integrable}
    The reader should compare our definition of almost complex structure of splitting type with the one of complex structure of splitting type given in \cite{Kas13}, see also \cite{AOUV17}. In addition to the assumptions given in Definition \ref{def:splitting}, when dealing with complex structure of splitting type one assumes that $H$ is nilpotent, that $J_H$ is integrable, and that the action of $\eta$ on the lie algebra of $H$ is semi-simple. This guarantees that the resulting structure on $G$ is integrable and that $G$ is solvable. The original definition also requires some special cohomological properties on $H$. We drop the assumption since in this paper we are not concerned with the cohomology of almost complex manifolds.
\end{remark}

\begin{lemma}\label{lemma:holomorphic:sum}
    Let $G = \C^k \ltimes_\eta H$ be an almost complex Lie group endowed with a \textbf{complex} structure of splitting type. Then the standard left-invariant $\su (m)$-structure on $G$ is holomorphic.
\end{lemma}
\begin{proof}
    The basis of $(1,0)$-forms for $J_H$ defines a left-invariant $(m-k,0)$-form
    \[
    \hat{\Omega} \coloneqq \psi^{k+1 \ldots m}
    \]
    on $H$. By \cite{BDV09}, if $H$ is nilpotent and $J_H$ is integrable, then $\bar \partial \hat{\Omega} =0$. In particular, since the action of $\eta$ is holomorphic with respect to $J_H$ and it can be diagonalized over $\C$ by semi-simplicity, we have that
    \[
    \bar \partial (\eta^{-1} \hat \Omega) = 0,
    \]
    proving the lemma.
\end{proof}

Almost complex structures of splitting type on $G$ and the associated left-invariant $\su (m)$-structure naturally descend to invariant structures on the compact quotient $M \coloneqq \Gamma \backslash G$. By assumption, the lattice $\Gamma$ splits as $\Gamma = \Gamma_{\C^k} \ltimes_\eta \Gamma_H$. Hence, there is a well-defined projection
\begin{equation}\label{eq:fibration}
\pi \colon M \longrightarrow T^{2k} \coloneqq \Gamma_{\C^k} \backslash \C^k,
\end{equation}
that gives to $M$ the structure of a fiber bundle over $T^{2k}$, with fiber diffeomorphic to $\Sigma \coloneqq \Gamma_H \backslash H$. By the expression of \eqref{eq:std:basis}, the fibration is pseudoholomorphic with respect to the almost complex structure of splitting type on $M$ and the complex structure on $T^{2k}$. 

\begin{example}
    Tori, as well as the direct products of  by $\C^k$ by almost complex Lie groups endowed with a lattice, trivially admit almost complex structures of splitting type. The Iwasawa manifold and the Nakamura manifold, which are two widely studied solvmanifolds endowed with complex and almost complex structures, both have many almost complex structures of splitting type, and provide non-trivial examples. We refer to Section \ref{sec:examples} for more non-trivial examples. 
\end{example}

Finally, due to the importance it plays in Theorem \ref{thm:family:0}, we denote by $\mathcal{{F}}$ the space of smooth functions on $M$ that are constant along the fibers of its fibration structure. In particular, if we denote by $\xi_j$, $j =1,\ldots,m$, a basis of $(1,0)$-vector fields dual to \eqref{eq:std:basis}, then we have that
\[
\mathcal{F} = C^\infty (M) \cap \ker \xi_{k+1} \cap \ldots \cap \ker \xi_m.
\]
The inclusion of $\mathcal{F}$ into the kernel of each operator $\xi_j$, $j=k+1, \ldots, m$, is clear since functions in $\mathcal{F}$ are constant on the fibers. Denote by $\partial_\Sigma$ the operator on the fibers of $\pi$ such that
\[
\partial f = \sum_{j=1}^k \partial_{z_j} (f) dz^j + \partial_\Sigma f,
\]
and let $f \in C^\infty (M) \cap \ker \xi_{k+1} \cap \ldots \cap \ker \xi_m$. This implies that $\partial_\Sigma f =0$. The operator $\bar \partial_\Sigma \partial_\Sigma$ is elliptic on $\Sigma$, which is compact. Hence, by the maximum principle, the function $f$ is constant on each fiber, proving the opposite inclusion.

\section{Kodaira dimension of almost complex structures of splitting type}\label{sec:fibrations}

Let $G$ be an almost complex Lie group of splitting type and let $J_0$ be an almost complex structure of splitting type on it. If the associated $\su (m)$-structure is pseudoholomorphic, then $K_{J_0}$ admits an invariant pseudoholomorphic trivialization that descends to the quotient $M = \Gamma \backslash G$. In particular, the structure $J_0$ induces on $M$ an invariant structure $J$ that has Kodaira dimension $0$.

In the next theorem, we describe a construction that allows to produce non-invariant structures with Kodaira dimension $0$.

\begin{theorem}\label{thm:family:0}
    Let $G = \C^k \ltimes_\eta H$ be a $2m$-dimensional Lie group endowed with an almost complex structure of splitting type whose standard $\su (m)$-structure is pseudoholomorphic, and let $M = \Gamma \backslash G$ be its compact quotient. Then there exists a family of non-invariant almost complex structures $J_F$ on $M$ with $\kappa_{J_F} =0$, parametrized by $F \in (\mathcal{F})^{m-k}$. 
\end{theorem}

\begin{proof}
    We first prove the theorem for the case $k=1$. Let $\phi^j$, $j=1, \ldots, m$, be a basis of invariant $(1,0)$-forms for the standard almost complex structure $J_0$ on $G$, with differentials
    \begin{equation}\label{eq:differentials}
    d \phi^j = \sum_{p <q } ( C^j_{p q} \phi^{p q} +C^j_{\bar p \bar q} \phi^{\bar p \bar q}) + \sum_{p,q} C^j_{p \bar q} \phi^{p \bar q}, \quad j =1, \ldots, m,
    \end{equation}
    with $C^j_{pq}$, $C^j_{\bar p \bar q}$ and $C^j_{p \bar q} \in \C$. Since $G$ is of splitting type, we can suppose $d \phi^1 =0$. Since the invariant $\su (m)$-structure associated to $J_0$ is pseudoholomorphic, we also have that 
    \[
    0 = \bar \partial (\phi^{1 \ldots m}) = \sum_j (-1)^{j+1} \bar \partial \phi^j \wedge \phi^{1 \ldots \hat{j} \ldots m} = - (\sum_{j,q} C^j_{j \bar q} \phi^{\bar l} ) \wedge \phi^{1 \ldots m}.
    \]
    In particular, this gives the relations
    \begin{equation}\label{eq:holomorphic}
    \sum_j C^j_{j \bar q} = 0, \quad \text{for } q =1, \ldots, m.
    \end{equation}
    Let $f_j \in \mathcal{F}$, $j = 2, \ldots, m$, be non-constant functions, and consider the non-invariant almost complex structure $J_F$ on $M$ defined by the $(1,0)$-forms
    \[
    \omega ^1 \coloneqq \phi^1 \quad \text{and} \quad \omega^j \coloneqq \phi^j + f_j \phi^{\bar 1} \quad \text{for } j = 2, \ldots, m,
    \]
    with $ F = (f_2, \ldots, f_m ) \in \mathcal{F}^{m-1}$. The inverse of the transformation is
    \[
    \phi^1 = \omega^1, \quad \phi^j = \omega^j - f_j \omega^{\bar 1} \quad \text{for } j = 2, \ldots, m.
    \]
    Moreover, if we denote by $\xi_j$, $j= 1, \ldots, m$, the $(1,0)$-vector fields dual to the $\phi^j$, then the $(1,0)$-vector fields dual to the $\omega^j$ are
    \[
    \psi_1 = \xi_1 - \sum_{j \ge 2} \bar{f}_j \xi_{\bar j} \quad \text{and} \quad \psi_j = \xi_j, \quad \text{for } j = 2, \ldots, m.
    \]
    To prove the theorem, we show that, for our choice of the functions $f_j$, we have that 
    \[
    \bar \partial_{F} (\omega^{1 \ldots m}) = \sum_{j} (-1)^{j+1} \, \bar \partial_F \omega^j \wedge \omega^{1 \ldots \hat{j} \ldots m} = 0,
    \]
    where $\bar \partial_F$ denotes the operator $\bar \partial$ associated to the structure defined by the $\omega^j$. To measure the difference of the action of $\bar \partial_F $ and $\bar \partial$, it is useful to define
    \[
    \bar \partial_0 \phi^j := \sum_{p,q} C^j_{p \bar q} \omega^{p \bar q},
    \]
    obtained by formally replacing $\phi^j$ by $\omega^j$ in the expression of $\bar \partial \phi^j$ given in \eqref{eq:differentials}. In particular, it follows from \eqref{eq:holomorphic} that
    \begin{equation}\label{eq:dbar0}
    \sum_j (-1)^{j+1} \, \bar \partial_0 \phi^j \wedge \omega^{1 \ldots \hat{j} \ldots m} = 0.
    \end{equation}
    Then, we have that
    \begin{equation}\label{eq:df:omega}
    \bar \partial_F \omega^1 = 0 \quad \text{and} \quad \bar \partial_F \omega ^j = \bar \partial_0 \phi^j + \epsilon_j + \partial_F f_j \wedge \omega^{\bar 1}, \quad \text{for } j = 2, \ldots, m,
    \end{equation}
    where $\epsilon_j$ measures the difference between $\bar \partial_F$ and $\bar \partial_0$. So we have that
    \begin{align*}
    \bar \partial_{F} (\omega^{1 \ldots m}) &= \sum_{j \ge 2} (-1)^{j+1} \, \bar \partial_F \omega^j \wedge \omega^{1 \ldots \hat{j} \ldots m}\\
    &= \sum_{j \ge 2} (-1)^{j+1} \, (\bar \partial_0 \phi^j + \epsilon_j + \partial_F f_j \wedge \omega^{\bar 1}) \wedge \omega^{1 \ldots \hat{j} \ldots m} \\
    &= \sum_{j \ge 2} (-1)^{j+1} \, \epsilon_j \wedge \omega^{1 \ldots \hat{j} \ldots m},
    \end{align*}
    where in the second equality we used \eqref{eq:df:omega}, while in the last equality we used \eqref{eq:dbar0} and the fact that $f_j \in \mathcal{F}$ gives
    \begin{equation}\label{eq:differential:difference}
    \partial_F f_j \wedge \omega^{1 \ldots \hat{j} \ldots m} = \psi_1 (f_j) \, \omega^1 \wedge \omega^{1 \ldots \hat{j} \ldots m} = 0.
    \end{equation}
    To obtain the explicit expression of the $\epsilon_j$, we write \eqref{eq:differentials} in terms of the $\omega^j$ and take the $(1,1)$-bidegree part. First, we have that
    \[
    \sum_{p <q } C^j_{pq} \, \phi^{pq} = \sum_{q > 1} C^j_{1q} \, \omega^1 \wedge (\omega^q- f_q \omega^{\bar 1}) + \sum_{1 < p <q} C^j_{pq} (\omega^p- f_p \omega^{\bar 1}) \wedge (\omega^q- f_q \omega^{\bar 1}).
    \]
    Taking bidegree $(1,1)$, the remaining terms are
    \begin{equation}\label{eq:extra:20}
    \sum_{q > 1} -f_q \, C^j_{1q} \, \omega^{1 \bar 1} + \sum_{1 < p <q} C^j_{pq} \, \omega^{\bar 1} \wedge (f_q \omega^p- f_p \omega^q).
    \end{equation}
    Similarly, when developing $\sum_{p < q} C^j_{\bar p \bar q} \, \phi^{\bar p \bar q}$, we get the extra terms
    \begin{equation}\label{eq:extra:02}
    \sum_{q > 1} \bar{f}_q \, C^j_{\bar 1 \bar q} \, \omega^{1 \bar 1} + \sum_{1 < p <q} C^j_{\bar p \bar q} \, \omega^1 \wedge (\bar{f}_q \omega^{\bar p}- \bar{f}_p \omega^{\bar q}).
    \end{equation}
    Finally, we have that 
    \begin{align*}
    \sum_{p,q} C^j_{p \bar q} \, \phi^{p \bar q} &= \sum_q (C^j_{1 \bar q} \, \omega^1 \wedge (\omega^{\bar q} - \bar{f}_q \omega^1) + C^j_{p \bar 1} (\omega^p -f_p \omega^{\bar 1}) \wedge \omega^{\bar 1})\\ 
    &+ \sum_{p, q \neq 1} C^j_{p \bar q} (\omega^p -f_p \omega^{\bar 1}) \wedge ( \omega^{\bar q} - \bar{f}_q \omega^1 ).
    \end{align*}
    Taking bidegree $(1,1)$ and using \eqref{eq:extra:20} and \eqref{eq:extra:02}, we have that
    \begin{align*}
    \bar \partial_F \phi^j &= \sum_{p,q} C^j_{p \bar q} \, \omega^{p \bar q} + \sum_{q > 1} (\bar{f}_q \, C^j_{\bar 1 \bar q} - f_q \, C^j_{1q}) \omega^{1 \bar 1}\\ 
    &+ \sum_{1 < p <q} C^j_{pq} \, \omega^{\bar 1} \wedge (f_q \omega^p- f_p \omega^q) + \sum_{1 < p <q} C^j_{\bar p \bar q} \, \omega^1 \wedge (\bar{f}_q \omega^{\bar p}- \bar{f}_p \omega^{\bar q}).
    \end{align*}
    By definition, the first summand is $\bar \partial_0 \phi^j$, while the second and fourth summand are irrelevant when wedging with $\omega^{1 \ldots \hat{j} \ldots m}$, since they all contain the term $\omega^1$. This leaves us with
    \begin{align*}
    \bar \partial_F (\omega^{1 \ldots m} ) &= \sum_{j \ge 2} (-1)^{j+1} \, \epsilon_j \wedge \omega^{1 \ldots \hat{j} \ldots m} \\
    &= \sum_{j \ge 2} (-1)^{j+1} \sum_{1 < p <q} C^j_{pq} \, \omega^{\bar 1} \wedge (f_q \omega^p- f_p \omega^q) \wedge \omega^{1 \ldots \hat{j} \ldots m}\\
    &= \sum_{j \ge 2} ( \sum_{ q > j} C^j_{jq} f_q - \sum_{1< p<q} C^j_{pj} f_p ) \omega^{\bar 1 1 \ldots m}.
    \end{align*}
    By rearranging the indices in a suitable way, we can write 
    \[
    \sum_{j \ge 2} ( \sum_{ q > j} C^j_{jq} f_q - \sum_{1 < q<j} C^j_{qj} f_q )  = \sum_{j \ge 2} A_j f_j,
    \]
    where $A_j \in \C$ has the expression
    \[
    A_j = \sum_{2 \le q < j} C^q_{qj} - \sum_{j < q \le m } C^q_{jq}, \quad j = 2, \ldots, m.
    \]
    To conclude the proof, we show that $A_j = 0$ for all $j$. Using the duality pairing between forms and vector fields on $\g$, and the fact that the differential of invariant forms is, up to sign, the dual of the Lie bracket, we have that
    \begin{align*}
    A_j &= \sum_{2 \le q < j } d \phi^q (\xi_q, \xi_j ) - \sum_{j < q \le m} d\phi^q (\xi_j, \xi_q) \\
    &= \sum_{q = 1}^m d\phi^q (\xi_q, \xi_j) = - \sum_{q=1}^m \phi^q ([\xi_q, \xi_j]) = \tr (\ad_{\phi^j}),
    \end{align*}
    where in the second equality we used that $[\xi_j, \xi_j ] = 0$ and that $[\xi_1, \xi_q] =0$ for all $q=1,\ldots,m$. Since $G$ admits a lattice, it is unimodular and the trace of the adjoint representation vanishes, giving $A_j=0$ for all $j =1,\ldots, m$.

    The proof for the case $k >1$ follows similarly, choosing $f_j =0$ for $j = 1, \ldots, k$. The only difference is in \eqref{eq:differential:difference}, that becomes
    \[
    \partial_F f_j \wedge \omega^{1 \ldots \hat{j} \ldots m} = \sum_p \psi_p (f_j) \, \omega^p \wedge \omega^{1 \ldots \hat{j} \ldots m} = (-1)^j \psi_j (f_j) \omega^{1 \ldots m}.
    \]
    The function $\psi_j(f_j)$ vanishes for every $j$ since we are taking $f_j =0$ when $j \le k$, and $f_j \in \mathcal{F}$ when $j \ge k+1$, proving the theorem for $k > 1$.
\end{proof}

In the second result, we describe a construction that allows to produce non-invariant structures with Kodaira dimension $-\infty$.

\begin{theorem}\label{thm:family:infty}
    Let $G = \C^k \ltimes_\eta H$ be a $2m$-dimensional Lie group endowed with an almost complex structure of splitting type whose standard $\su (m)$-structure is pseudoholomorphic, and let $M = \Gamma \backslash G$ be its compact quotient. Then there exists a family of non-invariant almost complex structures $J_F$ on $M$ with $\kappa_{J_F} = - \infty$, parametrized by $F \in C^\infty(M)^k$. 
\end{theorem}

\begin{proof}
    Let $\phi^j$, $j=1, \ldots, m$, be a basis of invariant $(1,0)$-forms for an almost complex structure of splitting type $J_0$ on $G$, with differentials given as \ref{eq:differentials}. Since $J_0$ is of splitting type, we can suppose $d \phi^j =0$ for $j = 1, \ldots ,k$. Since the invariant $\su (m)$-structure associated to $J_0$ is pseudoholomorphic, we also have that \ref{eq:holomorphic} holds.
    
    Let $f_j \in C^\infty (M)$, $j = 1, \ldots, k$, be non-constant functions, and consider the non-invariant almost complex structure $J_F$ on $M$ defined by the $(1,0)$-forms
    \[
    \omega ^j \coloneqq \phi^j, \quad \text{for $j =1, \ldots, m-1,$} \quad \text{and} \quad \omega^m \coloneqq \phi^m + \sum_{j = 1}^k f_j \phi^{\bar j},
    \]
    with $ F = (f_1, \ldots, f_k ) \in C^\infty (M)^k$. The inverse of the transformation is
    \[
    \phi^j = \omega^j, \quad \text{for $j =1, \ldots, m-1,$} \quad \text{and} \quad \phi^m = \omega^m - \sum_{j=1}^k f_j \omega^{\bar j}.
    \]
    Moreover, if we denote by $\xi_j$, $j= 1, \ldots, m$, the vector fields dual to the $\phi^j$, then the vector fields dual to the $\omega^j$ are
    \[
    \psi_j = \xi_j - \bar{f}_j \xi_{\bar m}, \quad \text{for $j =1, \ldots, m-1,$} \quad \text{and} \quad \psi_m = \xi_m.
    \]
    To prove the theorem, we start by showing that, for a suitable choice of the functions $f_j$, we have that the equation
    \[
    \bar \partial_{F} ( g \, \omega^{1 \ldots m}) = 0,
    \]
    has no solution $g \in C^\infty (M)$, where $\bar \partial_F$ denotes the operator $\bar \partial$ associated to the structure defined by the $\omega^j$. With the same notation of the proof of Theorem \ref{thm:family:0}, we have that
    \begin{align*}
    &\bar \partial_F \omega^j =  \partial_0 \phi^j + \epsilon_j, \quad \text{for } j = 1, \ldots, m-1, \quad \text{and} \\
    &\partial_F \omega ^m = \bar \partial_0 \phi^m + \epsilon_m + \sum_{j = 1}^k \partial_F f_j \wedge \omega^{\bar j} + \sum_{j=1}^k f_j \partial_F \omega^{\bar j}.
    \end{align*}
    Since $d \phi^j = d \omega^j = 0$ for $j = 1, \ldots, k$, the last summand in the differential $\partial_F \omega ^m$ is zero, so that
    \begin{align*}
    \bar \partial_{F} (\omega^{1 \ldots m}) &= \sum_{j = 1}^m (-1)^{j+1} \, \bar \partial_F \omega^j \wedge \omega^{1 \ldots \hat{j} \ldots m}\\
    &= \sum_{j = 1}^m  (-1)^{j+1} \, (\bar \partial_0 \phi^j + \epsilon_j) \wedge \omega^{1 \ldots \hat{j} \ldots m} + (-1)^{m+1} \sum_{j =1}^k \partial_F f_j \wedge \omega^{\bar j1 \ldots m-1} \\
    &= \sum_{j = 1}^m (-1)^{j+1} \, \epsilon_j \wedge \omega^{1 \ldots \hat{j} \ldots m} - \sum_{j = 1}^k \psi_m (f_j) \omega^{ \bar j} \wedge \omega^{1 \ldots m}
    \end{align*}
    where in the third equality we used \eqref{eq:dbar0}. To obtain the explicit expression of the $\epsilon_j$, we proceed as in the proof of Theorem \ref{thm:family:0}. The contribution of the terms $C^j_{p q} \, \phi^{pq}$ to the differential $\bar \partial_F \omega^j$ is 
    \[
    - \sum_{p < m} \sum_{q = 1}^k C^j_{pm} f_q \, \omega^{p \bar q}.
    \]
    Similarly, the terms $C^j_{\bar p \bar q} \, \phi^{\bar p \bar q}$ contribute as
    \[
    \sum_{p < m} \sum_{q = 1}^k C^j_{\bar p \bar m} \bar f_q \, \omega^{q \bar p}.
    \]
    Finally, the last contribution is that of $C^j_{m \bar m} \, \phi^{ m \bar m}$, that gives
    \[
    - \sum_{p,q= 1}^k C^j_{m \bar m} f_q \bar f_p \omega^{p \bar q}
    \]
    With a slight abuse of notation, and to write the sums in a more compact way, we set $f_j = 0$ for $j\ge k+1$. This leaves us with
    \begin{align*}
    \bar \partial_F (\omega^{1 \ldots m} ) &= \sum_{j = 1}^m (-1)^{j+1} \, (\sum_{p,q = 1}^{m-1} (-C^j_{pm} f_q + C^j_{\bar q \bar m} \bar f_p - C^j_{m \bar m} f_q \bar f_p) \omega^{p \bar q} ) \wedge \omega^{1 \ldots \hat{j} \ldots m}\\
    &- \sum_{j = 1}^k \psi_m (f_j) \omega^{ \bar j} \wedge \omega^{1 \ldots m} \\
    &= - \sum_{j = 1}^m  (\sum_{q = 1}^{m-1} (-C^j_{jm} f_q + C^j_{\bar q \bar m} \bar f_j - C^j_{m \bar m} f_q \bar f_j) \omega^{\bar q} ) \wedge \omega^{1 \ldots m} \\
    &- \sum_{j = 1}^k \psi_m (f_j) \omega^{ \bar j} \wedge \omega^{1 \ldots m} \\
    &= \sum_{q = 1}^m H_q \, \omega^{\bar q} \wedge \omega^{1 \ldots m}
    \end{align*}
    where in the second equality we used that the only summands that matter are those for which $p = j$, and where the $H_q \in C^\infty(M)$ are given by 
    \[
    H_q = \sum_{j=1}^m (C^j_{jm} f_q - C^j_{\bar q m} \bar f_j + C^j_{m \bar m} f_q \bar f_j) - \psi_m (f_q) \quad \text{for $q=1, \ldots, m-1$,}
    \]
    and $H_m = -\psi_m (f_m )$. The first summand vanishes since $\sum_{j=1}^m C^j_{jm}$ is the trace of the adjoint of $\phi^m$ and $G$ is unimodular. The second and third summands vanish since either $j \le k$, so that $C^j_{\bar q m} = 0$, or $j \ge k+1$, so that $f_j =0$. Recalling that we set $f_j =0$ for $j \ge k+1$, we conclude that
    \[
    H_q = - \psi_m (f_q), \quad \text{for $q = 1, \ldots, k$} \quad \text{and} \quad H_q = 0, \quad \text{for $q = k+1, \ldots, m$.}
    \]
    The equation $\bar \partial_F (g \omega^{1 \ldots m})$ has a solution $g \in C^\infty(M)$ if and only if $g$ solves the system
    \begin{equation}\label{eq:system}
    \begin{cases}
        \psi_q (g) = g \psi_m (f_q), & \text{for } q \le k, \\
        \psi_q (g) = 0, & \text{for }q \ge k+1.
    \end{cases}
    \end{equation}
    To show that it has no solution, observe that, by the last equation, we have that $\psi_m (g) =0$, that implies $\psi_q (g) = \xi_q (g)$ for $ q = 1, \ldots, m$. Consider the fibration $M \rightarrow T^{2k}$ defined by the splitting structure on $G$. Then the operator 
    \[
    \Delta_{k+1, \ldots, m} \coloneqq \sum_{j= k+1}^m \xi_j \xi_{\bar j}
    \]
    is elliptic on the fibers and any solution $g$ of \eqref{eq:system} lies in its kernel. By the maximum principle, the function $g$ must be constant on the fibers, and it descends to a function on the base $T^{2k}$. By the equations of \eqref{eq:system} for $q \le k$, if we choose the $f_q$ non-constant on the fibers, then the system has no solution $g$. This implies that the canonical bundle has no pseudoholomorphic section.

    To show that $\kappa_{J_F} = -\infty$, consider the equation $\bar \partial_F( g (\omega^{1 \ldots m} )^{\otimes l} )= 0$, for $l \ge 2$. Then $g$ is a solution of the equation if and only if it solves the system
    \[
    \begin{cases}
        \psi_q (g) = l \, g \psi_m (f_q), & q \le k, \\
        \psi_q (g) = 0, & q \ge k+1.
    \end{cases}
    \]
    With the same reasoning as above, we conclude that there are no solutions and that the Kodaira dimension of $J_F$ is $-\infty$. 
\end{proof}

\begin{cor}\label{cor:noninv}
    Every almost complex Lie group of splitting type whose standard $\su (m)$-structure is pseudoholomorphic admits non-invariant almost complex structures of Kodaira dimension $0$ and of Kodaira dimension $-\infty$.
\end{cor}

\begin{remark}
By Remark \ref{rem:integrable} and Lemma \ref{lemma:holomorphic:sum}, the claims of Theorems \ref{thm:family:0} and \ref{thm:family:infty}, and of Corollary \ref{cor:noninv} hold for complex structures of splitting type. They also hold for compact quotients of complex parallelizable Lie groups \cite{Wan54}, for the $6$-dimensional solvmanifolds with holomorphically trivial canonical bundle classified by Fino, Otal and Ugarte \cite{FOU15}, and for the compact quotients of the $6$-dimensional unimodular Lie groups with holomorphically trivial canonical bundle described by Otal and Ugarte \cite{OU23}.
\end{remark}

\section{Examples and applications}\label{sec:examples}

In this section we present several examples of Lie groups endowed with almost complex structures of splitting type and of their compact quotients. All the examples come from manifolds that, more or less extensively, have already been considered in the literature, showing the abundance of this kind of structure.

\begin{example}[\textbf{Almost complex structures with Kodaira dimension $-\infty$ on tori}]
    Complex tori fall into the assumptions of Theorems \ref{thm:family:0} and \ref{thm:family:infty}. When applying Theorem \ref{thm:family:0}, we have that $k =m$ and  we are not actually building a new structure. However, as a consequence of Theorem \ref{thm:family:infty} we are able to produce a large family of almost complex structure on tori of arbitrary dimension with Kodaira dimension $-\infty$. We point out that, alternatively, a different family of almost complex tori with Kodaira dimension $-\infty$ can be obtained by taking product of the $4$-dimensional example built by Chen and Zhang, see Section 6.2 in \cite{CZ23}, with tori endowed with an integrable structure.
\end{example}

\begin{example}[\textbf{The Iwasawa and the Nakamura manifolds}] The Iwasawa manifold and the Nakamura manifold are classical examples of solvmanifolds obtained as the quotient of $3$-dimensional complex Lie groups by a lattice \cite{Nak75}. They admit the structure of $T^4$-bundle over $T^2$, resp. of $T^2$-bundle over $T^4$, induced by a complex structure of splitting type. In particular, they all admit non-invariant almost complex structures with Kodaira dimension $0$ or $-\infty$.
    
\end{example}

\begin{example}[\textbf{A $6$-dimensional solvmanifold without invariant complex structures}]
The construction of the manifold of this example is due to Fern\'andez, de Le\'on and Saralegui \cite{FLS96}. Consider the Lie group given by the semi-direct product $G = \R^2 \ltimes_\eta \R^4$. If we denote by $(t,x)$ coordinates on $\R^2$, then $\eta(t,x)$ acts on vectors in $\R^4$ as left-multiplication by the matrix
\[
\eta(t,x) = \begin{bmatrix}
    e^t & 0 & x e^t & 0 \\
        & e^{-t} & 0 & x e^{-t} \\
        &          & e^t & 0 \\
        &          &     & e^{-t}
\end{bmatrix}.
\]
If we denote by $(y_1, y_2, z_1, z_2)$ coordinates on $\R^4$, a basis of left-invariant $1$-forms of $G$ is given by 
\begin{align*}
&e^1 = dt, \quad e^2 = dx, \quad e^3 = e^{-t} dy^1 - x e^{-t} dz^1,\\
&e^4 = e^t dy^2 - x e^t dz^2, \quad e^5 = e^{-t} dz^1, \quad  e^6 = e^{t} dz^2,
\end{align*}
with differentials
\begin{align*}
&de^1 = de^2 =0, \quad de^3 = -e^{13} - e^{25},\\
&de^4 = e^{14} - e^{26}, \quad de^5 = -e^{15} \quad \text{and} \quad de^6 = e^{16}.
\end{align*}

The Lie group $G$ admits a lattice $\Gamma$ and the quotient manifold $M = \Gamma \backslash G$ has the structure of a $T^4$-bundle over $T^2$. A prominent feature of $M$ is that it does not admit integrable invariant complex structures, even though it admits several structures of almost complex structures of splitting type induced by those on $G$. We refer to \cite{FLS96} for the details. We consider two invariant almost complex structures $J_1$ and $J_2$ on $M$ induced by different choices of integrable almost complex structure on $\R^4$. Both define an almost complex structure of splitting type on $G$ whose associated $\su(3)$-structure is pseudoholomorphic, so that we can apply Theorems \ref{thm:family:0} and \ref{thm:family:infty} to obtain families of non-invariant structures. However, the pseudoholomorphic trivialization of $K_{J_1}$ is invariant, while that of $K_{J_2}$ is necessarily non-invariant.

Let $J_1$ be the invariant almost complex structure on $M$ given by the co-frame of $(1,0)$-forms
\[
\phi^1 = e^1 + i e^2, \quad \phi^2 = e^3 + i e^4 \quad \text{and} \quad \phi^3 = e^5 + i e^6,
\]
with differentials
\[
d \phi^1 =0, \quad d\phi^2 = - \frac{1}{2} (\phi^{1 \bar 2} + \phi^{\bar 1 \bar 2} ) \quad \text{and} \quad d \phi^3 = -\frac{1}{2} (\phi^{1 \bar 3} + \phi^{\bar 1 \bar 3}) + \frac{i}{2} (\phi^{12} - \phi^{\bar 1 2}).
\]
Thus, we have that $\bar \partial \phi^{123} = 0 $ and that the standard $ \su (3)$-structure is pseudoholomorphic. By Theorems \ref{thm:family:0} and \ref{thm:family:infty}, we can find families of non-invariant almost complex structures on $M$, built starting from $J_1$, with Kodaira dimension $0$ or $-\infty$.

Let $J_2$ be the invariant almost complex structure on $M$ given by the co-frame of $(1,0)$-forms
\[
\phi^1 = e^1 + i e^2, \quad \phi^2 = e^3 + i e^5 \quad \text{and} \quad \phi^3 = e^4 + i e^6,
\]
with differentials
\begin{align*}
&d \phi^1 =0, \quad d\phi^2 = - \frac{1}{4} (\phi^{12} + \phi^{1 \bar 2} + 3 \phi^{\bar 1 2} - \phi^{\bar 1 \bar 2} ) \quad \text{and} \\
&d \phi^3 = \frac{1}{4} (3 \phi^{13} -  \phi^{1 \bar 3} + \phi^{\bar 1 3} - \phi^{\bar 1 \bar 3}).
\end{align*}
In this case, we have that $\bar \partial \phi^{123} = - \frac{1}{2} \phi^{\bar 1 1 2 3}$ and the canonical bundle is not trivialized by invariant forms. However, we now show that it is trivialized by a non-invariant form. 

Let $f \in C^\infty (M)$ and consider the equation $\bar \partial (f \phi^{123} ) =0$. This is equivalent to the vanishing of the $(0,1)$-form $\alpha = \bar \partial f - \frac{1}{2}f \phi^{\bar 1}$. If we denote by  $\xi_j$ a basis of $(1,0)$-vector fields dual to the $\phi^j$, then $\alpha=0$ if and only if $f$ is a solution of the system
\[
\begin{cases}
    \xi_{\bar 1} (f) = \frac{1}{2}f, \\
    \xi_{\bar 2} (f) = 0, \\
    \xi_{\bar 3} (f) = 0.
\end{cases}
\]
The last two equations, together with the $T^4$-bundle structure over $T^2$ of $M$, and the ellipticity on the fibers of the operator $\xi_2 \xi_{\bar 2} + \xi_3 \xi_{\bar 3}$, allow to conclude that $f$ depends only on the variables $(t,x)$. By writing $f = u + iv$, where $u$ and $v$ are real function defined on the basis $T^2$, and using that $\xi_{\bar 1} = \frac{1}{2}(\partial_t + i \partial_x)$, we are left with the system of real equations
\begin{equation}\label{eq:coupled:system}
\begin{cases}
    u_t - v_x - u = 0, \\
    u_x + v_t - v = 0.
\end{cases}
\end{equation}
Taking derivative with respect to $t$ of the first equation and with respect to $x$ of the second one, then summing the two equations, we get
\[
u_{tt} + u_{xx} - u_t -v_x = 0,
\]
and, by substituting the value of $v_x$ obtained from the first equation of \eqref{eq:coupled:system}, we have the decoupled equation
\begin{equation}\label{eq:decoupled}
u_{tt} + u_{xx} - 2u_t + u = 0.
\end{equation}
To find a solution, let 
\[
u = \sum_{m,n \in \Z} u_{mn} \, e^{imt + inx}, \quad u_{mn} \in \C,
\]
be the Fourier series expansion of $u$. Equation \eqref{eq:decoupled} is equivalent to the vanishing of
\[
u_{mn} ( -m^2 -n^2 -2im +1),
\]
for all $m,n$, so that the only non-zero coefficients of the Fourier series are $u_{01}$ and $u_{0-1}$, giving a solution
\[
u = A e^{ix} + B e^{-ix}, \quad A,B \in C.
\]
By the first equation of \eqref{eq:coupled:system}, we know that $v_x = -u$ and, by integration, we get
\[
v = i A e^{ix} - i B e^{-ix} + c(t),
\]
where $c(t)$ is a function depending only on $t$. Finally, from the second equation of \eqref{eq:coupled:system}, we must have
\[
0 = u_x +v_t - v = \partial_t (c(t)) - c(t).
\]
Since solutions of $c'(t) = c(t)$ are not periodic in $t$, we have that $c(t) = 0$ and that $f = 2 B e^{-ix}$, $B \in \C \setminus \{ 0 \}$, is a never-vanishing solution of the system that provides a pseudoholomorphic trivialization of the canonical bundle.
\end{example}

\begin{example}[\textbf{An $8$-dimensional solvmanifold}]

The construction of the manifold of this example is due to Sferruzza and the second author \cite{SfeT24}. Consider the Lie group given by the semi-direct product $G = \C \ltimes_\eta N$, where $N$ is the $3$-dimensional complex Heisenberg group. If we denote by $(z_1, z_2, z_3)$ coordinates on $N$, then $\eta (z_0)$ acts on $ (z_1, z_2, z_3)$ as left-multiplication by the matrix
\[
\begin{bmatrix}
    e^{-z_0} & & \\
     & e^{z_0} & \\
     & & 1
\end{bmatrix}.
\]
The almost complex structure of splitting type on $G$ induced by $dz^j$, $j=0,\ldots, 3$, is integrable and it is given by the co-frame of left-invariant $(1,0)$-forms
\[
\phi^1 = dz^0, \quad \phi^2 = e^{z_0} dz^1, \quad \phi^3 = e^{-z_0} dz^2 \quad \text{and} \quad \phi^4 = dz^3 -\frac{1}{2}z_1 dz^2 + \frac{1}{2} z_2 dz^1,
\]
with differentials
\[
d\phi^1 =0, \quad d\phi^2 = \phi^{12}, \quad d\phi^3 = - \phi^{13} \quad \text{and} \quad d\phi^4 = - \phi^{23}.
\]
In particular, the associated $\su(4)$-structure is holomorphic. Since $G$ admits a lattice $\Gamma$, the holomorphic $\su(4)$-structure descends to the quotient $M = \Gamma \backslash G$. We refer to Section 6 in \cite{SfeT24} for the construction of the lattice. By Theorems \ref{thm:family:0} and \ref{thm:family:infty}, there exist families of non-invariant almost complex structures on $M$ with Kodaira dimension $0$ or $-\infty$.
    
\end{example}

\printbibliography

\end{document}